\documentclass[letterpaper, 10 pt, conference]{ieeeconf}  
\IEEEoverridecommandlockouts          
\usepackage{amsmath}
\usepackage{amssymb} 
\usepackage{dsfont}

\newtheorem{remark}{Remark}

\usepackage{xcolor}
\usepackage{graphicx}
\usepackage{algorithm}
\usepackage{algpseudocode}
\usepackage[hyphens]{url}
\usepackage[hidelinks]{hyperref}
\usepackage{float}
\usepackage{comment}
\usepackage{caption}
\usepackage{subcaption}
\algrenewcommand\algorithmicrequire{\textbf{Input:}}
\algrenewcommand\algorithmicensure{\textbf{Output:}}
\usepackage{cite}

\title{\LARGE \bf
Reinforcement learning based demand charge minimization using energy storage
}

\author{Lucas Weber, Ana Bu\v{s}i\'c, Jiamin Zhu
\thanks{L. Weber is with Inria and DI ENS, École Normale Supérieure, CNRS, PSL Research University, Paris, France. E-mail: \href{mailto:lucas.weber@inria.fr}{lucas.weber@inria.fr}}
\thanks{A. Bu\v{s}i\'c is with Inria and DI ENS, École Normale Supérieure, CNRS, PSL Research University, Paris, France}
\thanks{J. Zhu is with IFP Energies nouvelles, 1 et 4 avenue de Bois-Préau, 92852 Rueil-Malmaison, France}
}

\begin{document}

\maketitle
\thispagestyle{empty}
\pagestyle{empty}

\begin{abstract}
Utilities have introduced demand charges to encourage customers to reduce their demand peaks, since a high peak may cause very high costs for both the utility and the consumer. 
We herein study the bill minimization problem for customers equipped with an energy storage device and a self-owned renewable energy production.
A model-free reinforcement learning algorithm is carefully designed to reduce both the energy charge and the demand charge of the consumer. The proposed algorithm does not need forecasting models for the energy demand and the renewable energy production. The resulting controller can be used online, and progressively improved with newly gathered data. 
The algorithm is validated on real data from an office building of IFPEN Solaize site. 
Numerical results show that our algorithm can reduce electricity bills with both daily and monthly demand charges.

\end{abstract}

\section{INTRODUCTION}
In the spot electricity market, the price is determined by the marginal cost of the most expensive power plant in operation \cite{Basics_Electricity_Price_Formation}. High peak demand often requires the use of more expensive and more polluting power generators, such as coal and natural gas \cite{Relation_price_pollution}. It is therefore in the interests of energy suppliers and society to reduce the peak demand to avoid the use of most expensive generators. 
Mechanisms to motivate customers to change their consumption behaviors are referred to as demand-side management \cite{Demand_Side_Management}. 
A common one for smoothing demand curves is the implementation of time-of-use (TOU) pricing \cite{Time_of_Use_pricing} with on-peak and off-peak hours. 
US utilities have further introduced demand charges (DC) to encourage large-load customers (usually commercial or industrial ones) to flatten their load curves \cite{Demand_Charges}. The electricity bill takes the form of
\begin{equation} \label{eqn: bill}
    fees + C_{energy}^{BP} + C_{demand}^{BP},
\end{equation}
where fees represent fixed subscription fees, $BP$ is a billing period that may be one day or one month, $C_{energy}^{BP}$ is the 
energy charge representing the cost of the total amount of energy exchanged with the grid, 
and $ C_{demand}^{BP}$ is the DC that is proportional to the maximum electricity demand during $BP$.
While \textit{monthly demand charge} (MDC) aims to reduce the calls to high-cost power plants, the \textit{daily demand charge} (DDC) was designed to additionally "incentivize customers with a behind-the-meter energy storage system (ESS) to shave peak loads throughout the entire month, to capture demand charge reduction savings more effectively "\cite{PGOptionStorage}.

In this article, we aim to provide an online controller to minimize the bill of the form \eqref{eqn: bill} for consumers equipped with a battery and a self-owned solar power generation. 
The main difficulties come from (i) the volatility of the renewable generations and the uncertain future demands, and (ii) the non-smoothness of peak demand (due to the max term).
Note that (ii) can be avoided by smooth reformulations, but at the price of increasing the state 
dimension \cite{Augmented_state_space}.

A common way to address the peak minimization problem is to manage the stored energy by stochastic optimization or model predictive control, e.g. \cite{Radiation_forecasts,Augmented_state_space}. In such methods, future energy production and consumption must be forecasted, and their performance relies heavily on the quality of these forecasts. 
Unfortunately, it is difficult to build precise forecasting models for renewable power generation \cite{Energy_forecasts} and individual electricity demand \cite{Demand_forecasts}.
An alternative to avoid forecasting is the use of robust optimization \cite{xiang2015robust}, which emphasizes the worst-case performance under the uncertainty set, and hence may appear to be too conservative.
\cite{mo2021optimal} further proposed an online algorithm based on competitive ratio analysis. The proposed algorithm manages to discharge the energy storage to achieve the best competitive ratio.
However, the analysis was done with a linear battery model, and only the battery discharge is considered.

In this article, we propose to apply another model-free alternative: reinforcement learning (RL).
RL agents learn optimal control strategies from observations of the return of the environment.
RL has been successfully applied to minimize energy charges in various settings (see e.g. \cite{EC1, EC2, EC3}),
however, accounting for DC in RL remains challenging.
We design a Decomposition-based Policy Iteration (DPI) algorithm, which is applied off-line to obtain an optimal state-action value function ($Q$-function) on the training set. This $Q$-function can then be used as a heuristic by an on-line controller. The trained controllers can also be progressively improved (re-trained) with new collected data, and thus adapt to either new demand patterns or production changes (e.g. caused by climate changes).
In contrary to most existing solutions, we consider a nonlinear model of the battery, which makes the problem harder.
To be able to deal with the monthly charges, we split the training into separate days and suppose an empty battery at the end of each day. 
Our algorithm exploits this decomposition to accelerate traditional discrete Policy Iteration (PI) algorithms: an average of exact Dynamic Programming (DP) solutions (for separate training days) is used to initialize the PI algorithm, and a decomposition-based policy evaluation step is designed to accelerate the training.
We have adopted the tabular setting in this article, since it enjoys convergence guarantees, and is simple to implement without the need for fine hyperparameter tunings.

The remainder of this paper is organized as follows. Section \ref{section 3} provides background on RL. Section \ref{section 4} formalizes electricity bill minimization problems and presents our reinforcement learning approach. Section \ref{section 5} shows economic benefits of our method in several contexts. Section \ref{section 6} summarizes our results and discusses possible improvements.

\section{RL background}
\label{section 3}

\subsection{Reinforcement learning}
In RL, an agent (i.e. decision-maker) 
learns
an optimal behavior 
from the costs (or rewards) received from its interactions with an environment. 
Formally, we define a finite horizon Markov Decision Process (MDP) $(\mathcal{T},\mathcal{S},\mathcal{A},\mathcal{C},\mathcal{P})$, where $\mathcal{T}=\{0, 1, \dots, H-1\}$ is a finite time horizon,
$\mathcal{S}$ is the state space, $\mathcal{A}$ is the action space, $\mathcal{P}$ is the matrix of
transition probabilities of the environment, and $\mathcal{C}:\mathcal{T}\times\mathcal{S}\times\mathcal{A}\times\mathcal{S}\to \mathbb{R}$ is a cost function.
We denote $\mathcal{A}(s)$ the allowed actions in state $s$.
Action $a$ chosen by an agent is dictated by its policy. In this paper, we only consider deterministic policies. In state $s$ at time $\tau$, we write directly the action when following policy $\pi$ as $a = \pi(\tau,s)\in \mathcal{A}(s)$. 

The agent aims to find a policy $\pi^*$ minimizing the expectation of its finite horizon return:
\begin{equation*}
    \label{eq: expected return}
    \forall \tau\in\mathcal{T}, \quad \pi^*=\underset{\pi}{\arg \min}\; \mathbb{E}\left[ \sum_{k=\tau}^{H-1} \mathcal{C}_k(s_k, a_k, s_{k+1}) \right] 
\end{equation*}
where $\{s_k,a_k \}_{k=0,\cdots,H-1}\cup\{s_H\}$ is an agent's trajectory following policy $\pi$: $a_k = \pi(k, s_k)$ and $s_{k+1} \sim \mathcal{P}(s_k, a_k)$. 
For a policy $\pi$, we define the $Q$-function as
\begin{equation*}
\hspace*{-0.31cm} 
    Q^\pi(\tau,s,a)=\mathbb{E}\left[ \sum_{k=\tau}^{H-1} \mathcal{C}_k(s_k, a_k, s_{k+1}) \Big| s_\tau = s, a_\tau=a\right]
\end{equation*}
with $a_k = \pi(k,s_k) $.
For each time $\tau$ and state-action couple $(s,a)$, $Q^\pi(\tau,s,a)$ is the expected value of choosing action $a$ in state $s$ starting at time $\tau$ and following policy $\pi$ afterwards.
The optimal $Q$-function $Q^*$ is defined as:
\begin{equation*}
    \label{eq: Q^*}
    \forall (\tau,s,a) \in \mathcal{T}\times\mathcal{S}\times \mathcal{A}(s), \quad 
    Q^*(\tau,s,a) = \underset{\pi}{\min}\, Q^\pi(\tau,s,a),
\end{equation*}
and $Q^*$ satisfies the fixed point of Bellman optimality equation:
\begin{equation}
    \label{eq: Bellman optimality equation}
    Q^*(\tau, s,a)=\mathbb{E}_{s'}\left[ \mathcal{C}_\tau(s,a,s') + \underset{\tilde{a}\in \mathcal{A}(s)}{\min} Q^*(\tau', s', \tilde{a})\, \right]
\end{equation}
with $s' \sim \mathcal{P}(s,a)$.
Then, the optimal policy $\pi^*$ can be derived by selecting for any state $s$ and time $t$ the action $\pi^*(\tau, s)=\underset{a\in\mathcal{A}(s)}{\arg \min}\, Q^*(\tau, s, a)$.

To solve \eqref{eq: Bellman optimality equation}, two 
popular algorithms are Value Iteration (VI) and PI.
In this paper, we design a PI algorithm to solve the bill-minimization problem. The main advantage for using PI over VI is that it is empirically known to converge faster to an optimal policy.

The PI algorithms \cite[section 4.3]{Policy_Iteration} start with an initial policy $\pi_0$, and update it recursively by alternating the following evaluation and improvement steps.

\paragraph{Evaluation step} compute the fixed-policy $Q$-function $Q^{\pi_k}(\tau, s,a)$ that satisfies the fixed-point equation 
\begin{equation}
    \label{eq: Bellman equation with policy}
    Q^{\pi_k}(\tau, s,a)=\mathbb{E}_{s'}\left[ \mathcal{C}_\tau(s,a,s') + Q^{\pi_k}(\tau', s', a')\, \right]
\end{equation}
where $\tau'=\tau+1$, $Q^{\pi_k}(H, \cdot, \cdot)=0$, $s'\sim \mathcal{P}(s,a)$ and $a' = \pi_k(\tau', s')$. 
TD learning algorithms, such as SARSA, can be used to solve \eqref{eq: Bellman equation with policy}.

\paragraph{Policy improvement step} construct a new policy by
\begin{equation}
    \label{eq: policy improvement}
    \pi_{k+1}(\tau,s)=\underset{\tilde{a}\in \mathcal{A}(s)}{\arg \min}\, Q^{\pi_k}(\tau,s,\tilde{a})
\end{equation} for all time $\tau\in\mathcal{T}$ and state $s\in \mathcal{S}$.
The PI algorithm converges to an optimal policy \cite[section 4.3]{Policy_Iteration}.

\section{Application to bill minimization problems}
\label{section 4}
\subsection{Bill minimization problems}
\begin{figure}[H]
	\centering
  	\includegraphics[scale=0.4]{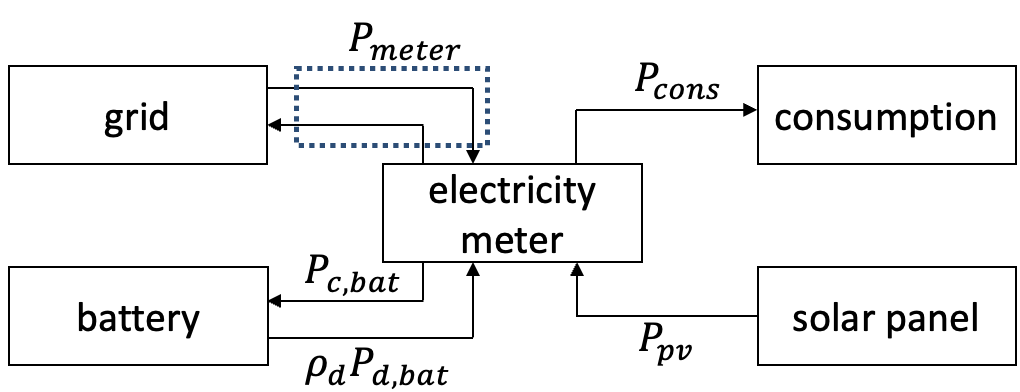}
   \caption{Sketch of the system}
    \label{fig: system}
\end{figure}

We consider a consumer equipped with a solar panel and a battery whose interactions with the grid are illustrated in Fig. \ref{fig: system}. The electricity meter centralizes the power flows. It supplies the consumption $P_{cons} [kW]$ to the user and   
 the charging power $P_{c,bat} [kW]$ to the battery and receives the solar power production $P_{pv} [kW]$ from the solar panel, the discharging power $\rho_d P_{P_{d,bat}} [kW]$ from the battery and $P_{meter}$ from the grid. $\rho_d [-]$ is the discharge efficiency of the battery: the electricity meter only receives a fraction $\rho_d$ of what the battery provides. Thus, the total balance power is given by
\begin{equation}
    \label{eq: Pmeter}
    P_{meter} = P_{cons} - P_{pv} + P_{c, bat} - \rho_d P_{d, bat}
\end{equation}
When $P_{meter}$ is positive, that amount of power $P_{meter}^+=\max(0, P_{meter})$ is bought from the grid, otherwise, the power $P_{meter}^-=-\min(0, P_{meter})$ is sold. 
$P_{meter}^+$ is indeed the net power demand.

Our objective is to minimize electricity bills of the form \eqref{eqn: bill}, 
where the energy charge and the DC are defined formally by
\begin{equation} \label{energy_charge}
    C_{energy}^{BP} = \int_{t\in BP} pp(t)P^+_\text{meter}(t) - sp(t)P^-_\text{meter}(t) \, dt
\end{equation}
\begin{equation} \label{demand_charge}
    C_{demand}^{BP} = \mu \; \underset{t\in BP}{\max} \left\{P^+_\text{meter}(t) \right\}.
\end{equation}
where $P_{meter}$ is defined by \eqref{eq: Pmeter},
$pp[\$ / kWh]$ and $sp [\$ / kWh]$ are the purchase and selling prices of energy, and $\mu [\$/kW]$ is a price per kilowatt. 
In this paper we consider a known $\mu$ (constant for simplicity) and a TOU tariff with fixed on-peak hours: $pp(\cdot)$ and $sp(\cdot)$ are known functions of time.

\medskip

In the rest of this paper, we are to study three bill minimization problems: 
\begin{itemize}
    \item daily energy minimization (DEM) : minimize only the energy charge $C^{BP}_{energy}$  \eqref{energy_charge} with the BP being a day, i.e., $BP=day$.
    \item DDC minimization (DDM) : minimize $C^{BP}_{energy} + C^{BP}_{demand}$ defined by \eqref{energy_charge}-\eqref{demand_charge} with $BP=day$.
    \item MDC minimization (MDM) : minimize $C^{BP}_{energy} + C^{BP}_{demand}$ defined by \eqref{energy_charge}-\eqref{demand_charge} with $BP=month$.
\end{itemize}
Note that the first one has been studied as an intermediate step towards DC minimization problems. Since subscription fees are fixed, they are excluded from the analysis. 

Problems (DEM), (DDM), and (MDM) are by nature multi-stage stochastic optimization problems, in which the matrix of transitions $\mathcal{P}$ is only partially known.
Indeed, we suppose the dynamic of the battery known, but not the evolution of the consumption and production. 
In next sections, we are to design a model-free RL algorithm to solve these problems. 
The proposed algorithm learns off-line and results in a (feedback-)policy that can be used on-line. 
Note that it can also be adaptive, in the sense that it can improve its policy gradually when new data of $P_{pv}$ and $P_{cons}$ are collected.

\subsection{Markov decision process}
\subsubsection{System model}
For our bill minimization problems, we consider a nonlinear space model, in which only a part of the model (battery dynamics) is known.
Recall that the transition properties for the production $P_{pv}$ and the consumption $P_{cons}$ are not known.  

Denoting by $x$ the state of charge (SOC) of the battery, the dynamics of the battery \cite{Dynamic_Battery} is given by
\begin{equation}
\label{eqn:x_dot}
\dot{x}(t)=I_{bat}(t) / Q_{bat}^{nominal}
\quad \left[s^{-1}\right]
\end{equation}
with $Q_{bat}^{nominal}[As]$ the nominal charge capacity of the battery and $I_{bat} [A]$ the current of the battery given by
\begin{equation*}
\label{eqn:Ibat}
\begin{cases}
\frac{-U_{ocv}(x)+\sqrt{U_{ocv}(x)^2+4R_{c,int}(x)P_{c,bat}(t)}}{2R_{c,int}(x)} & \text{(charge)}\\
\frac{-U_{ocv}(x)+\sqrt{U_{ocv}(x)^2-4R_{d,int}(x)P_{d,bat}(t)}}{2R_{d,int}(x)} & \text{(discharge)}
\end{cases}
\end{equation*}
The open circuit voltage $U_{ocv} [V]$ and the internal charging and discharging resistances $R_{c, int} [Ohm]$ and $R_{d,int} [Ohm]$ depend on $x$. 
The system controls are the charging and discharging powers of the battery $P_{c, bat} [kW]$ and $P_{d, bat}$ $[kW]$. 
The state of charge and the control of the system must satisfy the following constraints:
\begin{align}
\label{eqn: constraints 1}
& P_{c, bat} \in \left[0 , P_{c, bat}^{max}(x)\right], \quad P_{d, bat}  \in \left[0, P_{d, bat}^{max}(x)\right] \\
\label{eqn: constraints 3}
& x \in [0,1]
\end{align}
where the maximal charging and discharging powers $P_{c,bat}^{max} [kW]$ and $P_{d,bat}^{max}$ $[kW]$ depend on $x$.

\begin{remark}
 Such a model implies non-linear transfer efficiencies. In particular, the charging speed $\dot{x}$ is not a linear function of the (dis)charging power and the SOC. These realistic conditions may have an influence on the optimal strategies.
\end{remark}

\subsubsection{State and action spaces}
We consider the state space consisting of
\begin{itemize}
    \item the time step $\tau \in \mathcal{T}=\{0, 1, \dots, H-1\}$
    \item the state of charge $x \in \mathcal{S}o\mathcal{C} = [0, 1]$
    \item the highest demand peak registered so far during the current BP $p\in \mathcal{P} = [0, p_{max}]$
    \item the difference between the electricity demand and the solar production $\delta=P_{cons}-P_{pv} \in \Delta = [-\delta_{min}, \delta_{max}]$. 
\end{itemize}
and the action space consisting of $a=P_{c,bat}-P_{d,bat} \in \mathcal{A}$. 
By definition, we have $P_{c,bat} = a^{+}$ and $P_{d,bat} = a^{-}$.
As we have observed that small charging or discharging powers have higher battery efficiencies, we simplify the mixed state-control constraint \eqref{eqn: constraints 1} by a simple bound constraint $\mathcal{A}=[a_{min}, a_{max}]$ with 
$a_{max} = \min_{x\in[0,1]} P_{c,bat}^{max}(x)$ and $a_{min} = -\min_{x\in[0,1]}  P_{d,bat}^{max}(x)$.

\subsubsection{Cost function design}
To apply RL algorithms, immediate cost must be properly designed. 
In problems (DDM) and (MDM), the DC \eqref{demand_charge} constitutes a sparse signal that depends on the entire BP history.
Inspired by \cite{Augmented_state_space}, we design cost functions that only depends on the current state and action, as follows:
\begin{itemize}
    \item problem (DEM): a natural choice of the immediate cost function is the energy charge of one time step $\Delta t$, i.e.,
\begin{equation}
        \begin{aligned}
            \label{eq: cost function crit 1}
            \mathcal{C}^1_\tau(x, a)&= \int_{t=\tau\Delta t}^{(\tau+1)\Delta t} pp(\tau)P^+_{meter}(t, a) \, dt \\ & -\int_{t=\tau\Delta t}^{(\tau+1)\Delta t} sp(\tau)P^-_{meter}(t, a) \, dt
        \end{aligned}
    \end{equation}

    \item problem (DDM): we define the immediate cost as:
\begin{equation}
    \label{eq: cost function crit 2}
    \mathcal{C}^2_\tau(p, x, a)=\mathcal{C}^1_\tau(x, a) + (p'-p)^+\mu_d
\end{equation} 
with $\mu_d [\$/kW]$ is the DC rate for DDC, $p'$ the new maximal measured peak after the transition. 

    \item problem (MDM): the immediate cost is designed as :
\begin{equation}
    \label{eq: cost function crit 3}
    \mathcal{C}^3_\tau(p, x, a)=\mathcal{C}^1_\tau(x, a) + (p'-p)^+\mu_m
\end{equation} 
with $\mu_m [\$/kW]$ is the MDC rate.
\end{itemize}

\medskip

The second term in \eqref{eq: cost function crit 2}-\eqref{eq: cost function crit 3} are designed to overcome the noncumulative nature of demand peaks.
Since the real demand peak can be observed only at the end of the BP, we define the immediate cost of the DC by the increment of the peak $(p^{'}-p)^+$: an action leading to a higher peak has a larger cost. 
By doing so, at the end of the BP, the DC does correspond to the highest peak.

\subsection{Search for an optimal policy }
\subsubsection{Problem decomposition}
\label{problem decomposition}
For problems (DEM) and (DDM), a natural time horizon is one day. 
However, for the problem (MDM), the time horizon should be one month, which makes the problem much more difficult to solve. 
We propose therefore a problem decomposition by splitting the training into separate days: the agents are trained on a finite horizon $H$ corresponding to one day, starting at midnight, with an empty battery. 
Indeed, from a daily bill point of view, it is optimal for an agent to have an empty battery at the end of the day. If the battery was not empty, the remaining energy could have met the electricity demand or been sold. 

Note that this decomposition may result in suboptimal actions. If the battery sells electricity at the end of one day and buys it back at the beginning of the next, and the selling prices are lower than the purchase prices, it would have been more economical to store that energy. 
We will discuss the issue of recoupling later in section \ref{recoupling subsection}.

 With this problem decomposition for (MDM), the trainings for (DDM) and (MDM) are identical except for the price per kilowatt $\mu$.
 For both, the optimal $Q$-function is computed through DPI.
 However, the online applications differ: the way to reset the recorded peak $p$ in \eqref{eq: cost function crit 2}-\eqref{eq: cost function crit 3} is different for problems (DDM) and (MDM). For problem (DDM), $p$ is reset to zero at the beginning of each day, while for problem (MDM), $p$ is reset to zero only at the beginning of the month.
 
\medskip

\subsubsection{Decomposition-based PI algorithm}
\label{subsection: search for an optimal policy}
In this section, we design a Decomposition-based Policy Iteration algorithm.
The policy improvement step is classical, given by \eqref{eq: policy improvement}. 
However, the evaluation step is designed based on the decomposition described above.

Given a set of days $\mathcal{D}$ and $N=|\mathcal{D}|$ the number of days, at iteration $k$, the designed evaluation step is given by
\begin{equation} \label{eq:q_estimate}
    \hat{Q}^{\pi_k} = \mathbb{E}_\mathcal{D} \left[ Q_{d_i}^{\pi_k} \right]
\end{equation}
where $Q_{d_i}^{\pi_k}$ is defined as the $Q$-function for policy $\pi_k$ knowing the power consumption $P_{cons}$ and production $P_{pv}$ of day $d_i \in \mathcal{D}$, i.e.,
\begin{equation*}
    \label{eq: Q-function for day d}   Q_{d_i}^{\pi_k}=\mathbb{E}\left[ Q^{\pi_k} | P^{d_i}_{pv}, P^{d_i}_{cons} \right]
\end{equation*}
For a given policy $\pi_k$ and a given day $d_i$ where  $P^{d_i}_{pv}(\cdot)$ and $P^{d_i}_{cons}(\cdot)$ data are observed, $Q_{d_i}^{\pi_k}$ can be calculated by a classical DP
since the transitions are then deterministic.


Now, to validate the DPI algorithm, it remains to verify that, as $N \to \infty$, the estimate \eqref{eq:q_estimate} converges to the true $Q$-function $Q^\pi$. 
From the tower property, $\mathbb{E}_{d}\left[ Q^\pi_d \right]=Q^\pi$
and by construction, $\lim_{N\to\infty} \hat{Q}^{\pi} = \mathbb{E}_{d}\left[ Q^\pi_d \right]=Q^\pi$. 
\medskip

\subsubsection{Choice of initial policy} 
Recall that an PI algorithm begins with an initial policy $\pi_0$.
To speed up the training, we can provide an initial policy that may be close to the optimal policy $\pi^*$.  

In fact, given the set of training days $\mathcal{D}$ with $N$ the number of training days, the consumption $P_{cons}$ and the production $P_{pv}$ are known from the point of view of a single day, and  we can compute the optimal $Q$-function $Q_d^{\pi_d^*}$ for each of these days by classical DP-based methods, where $\pi_d^*$ is an optimal policy for day $d$.
Intuitively, an average of these optimal $Q$-functions $Q_d^{\pi_d^*}$ can be a good estimation of the real optimal $Q$-function $Q^*$.
We define 
$ \bar{Q}=\mathbb{E}_d\left[ Q_d^{\pi_d^*} \right]$.
This average can be estimated iteratively as follows: when state $s$ is visited at time step $\tau$ and action $a$ is chosen for day $d_n$,
\begin{equation} \label{q_average} 
    \begin{aligned}
    \bar{Q}_{n+1}(\tau, s, a) = & \left(1-\frac{1}{\text{nb visits of }(\tau, s, a)} \right)\bar{Q}_{n}(\tau, s, a) \\
    &+ \frac{1}{\text{nb visits of }(\tau, s, a)} Q_{d_n}^{\pi_{d_n}^*}(\tau, s, a)
    \end{aligned}
\end{equation}
with $\bar{Q}_{0}=0$.
Then, we can choose 
$$\pi_0(\tau, s)=\underset{\tilde{a}}{\arg \min}\, \bar{Q}_{N}(\tau, s, \tilde{a})
$$ 
as initialization for the PI algorithm.

\begin{remark}

The function $\bar{Q}$ is not necessarily a fixed-point of the Bellman optimality equation \eqref{eq: Bellman optimality equation}, since
\begin{equation*}
    \begin{aligned}
     \bar{Q}&(\tau,  s, a)  =  \mathbb{E}_{d} \left[ Q_d^{\pi_d^*}(\tau, s, a)\right] \\ &= \mathbb{E}_{d}\left[ \mathbb{E}_{s'}\left[ \mathcal{C}_{\tau}(s,a,s^\prime) + \underset{\tilde{a}}{\min}\, Q^{\pi_d^*}(\tau', s', \tilde{a}) \Big| P^d_{pv}, P^d_{cons}  \right]  \right] \\
     &\leq \mathbb{E}_{s'} \left[\mathcal{C}_{\tau}(s,a,s^\prime) + \underset{\tilde{a}}{\min}\,\bar{Q}(\tau^\prime,s',\tilde{a}) \right] 
     \end{aligned}
\end{equation*}
Therefore, the policy $\pi_0$ derived from $\bar{Q}$ is likely not optimal.
Nevertheless, we expect $\bar{Q}$ to be relatively close to the optimal $Q$-function $Q^*$. 
\end{remark}

\medskip

\subsubsection{Training and application}
Our DPI is trained off-line with a set of training days $\mathcal{D}$, and results in an estimate of the optimal $Q$-function $\hat{Q}^{*}$. 
Then, $\hat{Q}^{*}$ can be used online to choose the action $a_i$ given any state $x_i$ encountered at stage (or time step) $i \in \mathcal{T}$ by 
$$
a_i \in \underset{\tilde{a}\in \mathcal{A}}{\arg \min}\, \hat{Q}^{*}(\tau_i, s_i, \tilde{a})
$$
Again, we note that the agents can be frequently retrained to be able to adapt to the change of the environment. To do so, we only need to 
add recent days with newly collected data to the set $\mathcal{D}$ or replace several old training days.

During the training, at each time step, DP is performed for all possible values of $p$. 
During the test (or online application) phase, the RL agent can meet states that were never visited during the training.
We will compare to ways to handle this situation
\begin{itemize}
    \item Lazy controller: doing nothing by choosing $a=0$.
    \item Heuristic controller: the battery is charged when $P_{pv} \geq P_{cons}$. If it is full, the excess power is sold. The battery is discharged when $P_{pv} < P_{cons}$, but only during on-peak hours. At the end of the day, the battery is emptied.
\end{itemize}

\section{Numerical experiments and results}
\label{section 5}
We test our approach with a building of IFPEN's Solaize site equipped with a solar panel and a battery.
We have access to two years data of PV production and power demand, with a sample time of 10 minutes.
Since the power demands vastly differ between working and non-working days, we consider them separately. 
In this paper, we present results on working days, where a much greater bill reduction is possible compared to non-working days. For problem (MDM), we consider BP of $20$ working days. We provide experimental code here: \url{https://github.com/luweber21/dc-minimization}.

\subsection{Numerical settings}
The battery used for our numerical experiment has a storage capacity of $27.3$ kWh, costing roughly $\$4000$ given the current price per kWh. For the battery model \eqref{eqn:x_dot}, experimental data of $U_{ocv}$, $R_{c,int}$, $R_{d,int}$, $P_{c,bat,max}$ and $P_{d,bat,max}$ are used. $Q_{bat}^{nominal}$ is calculated from these measurements.
As mentioned, we use a fixed TOU tariff: the on-peak hours are from 6am to 9am and 6pm to 9pm with the price of $0.1936 \$/kWh$, while the off-peak price is set to be $0.1330 \$/kWh$, and the selling price is $sp = 0.098 \$/kWh$.
The DDC rate (resp. MDC rate) is $\mu_d=0.5 \$/kW$ (resp. $\mu_m=20\mu_d=10 \$/kW$).

We choose for the time space $\mathcal{T}=[0,143]$ which corresponds a horizon of one day with a time step of $\Delta t=10$ minutes. So defined, the time space $\mathcal{T}$ coincides with the sampling times of the average electricity demand $P_{cons}$ and average electricity production $P_{pv}$ measurements. In other words, time step $\tau$ corresponds to real time $t=\tau\Delta t$. 
Other spaces $\mathcal{S}o\mathcal{C}=[0,1]$, $\Delta=[-60,60]$, $\mathcal{P}=[0,100]$, and $\mathcal{A}=[-20,20]$ are discretized with a step of $0.01$, $2$, $1$, and $1$ respectively.
We denote these discretized spaces by $\mathcal{S}o\mathcal{C}_d$, $\Delta_d$, $\mathcal{P}_d$ and $\mathcal{A}_d$.

\medskip

In the next sections, we present results of agents that are trained over $N_{train}=300$ working days 
and tested over $N_{test}=100$ working days. 
The bill reduction is defined as the difference between the bill in the no-battery case and the optimized one.

\subsection{Better choice of state-spaces from solving problem (DEM) }
In this section, we study the problem (DEM) to determine the best state space for other experiments. 
Four agents $A_1-A_4$ with different state spaces are trained:
\begin{itemize}
    \item agent $A_1$ has state space $S_1=\mathcal{T} \times \mathcal{S}o\mathcal{C}_d$
    \item agent $A_2$ has state space $S_2=\mathcal{T} \times \mathcal{S}o\mathcal{C}_d \times \{-1, 0, 1\}$. The set $ \{-1, 0, 1\}$ is a simplification of $\Delta_d$, which takes the value of $1$ if $P_{pv} > P_{cons}$, $-1$ if $P_{pv} < P_{cons}$, and $0$ otherwise.
    \item agent $A_3$ has state space $S_3=\mathcal{T}
    \times \mathcal{S}o\mathcal{C}_d \times \Delta_d$, and is completed with the lazy controller during the tests.
    \item agent $A_4$: trained as $A_3$, but use the heuristic controller to complete the learned policy during the tests.
\end{itemize}

\begin{figure}
    \centering
    \includegraphics[width=0.8\columnwidth]{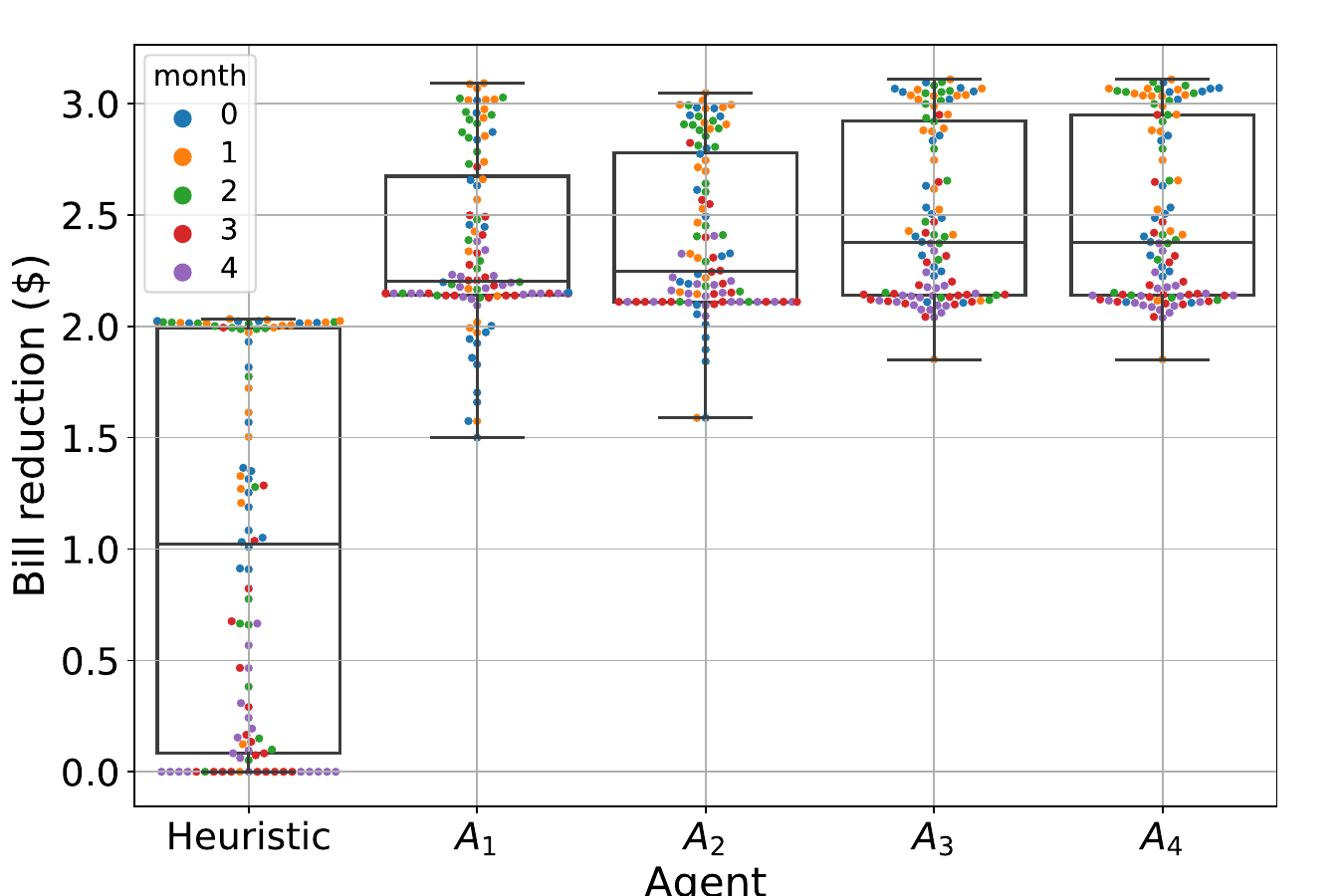}
    \caption{Swarm plot representing the daily bill reductions on the test set when using our agents versus no battery. 
    Each color represents one month.}
	\label{fig:bill minimization}
\end{figure}

For a discrete state space $\mathcal{S}$, the computational complexity of a DPI loop is $O\left(N_{train} |\mathcal{S}|\right)$.
In Fig. \ref{fig:bill minimization}, we show for each day of our test set the bill reduction obtained with the heuristic controller only and agents $A_1$-$A_4$.  
We observe that the average profits increase non-trivially with the complexity of the state spaces. 
Since state space $S_3$ gives the best results, we will use it as a basis for all other experiments. 

Comparing agent $A_3$ and agent $A_4$, the addition of the heuristic brings negligible additional benefits. It indicates that few states have not been visited during the training.

\subsection{Results of problem (DDM)}
\label{subsection: daily DC}
\begin{figure}
    \centering
    \includegraphics[width=0.8\columnwidth]{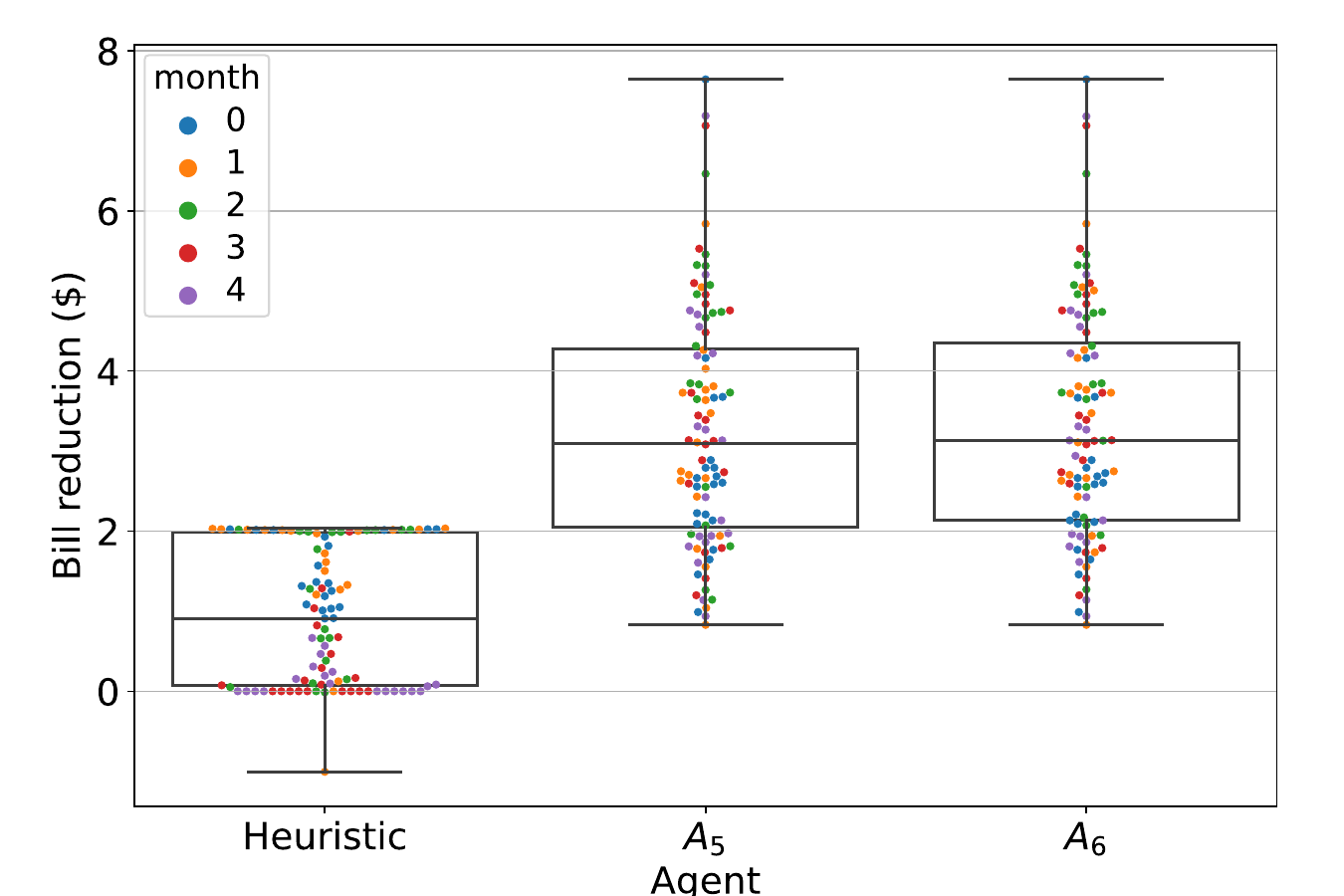}
    \caption{Swarm plot representing the bill reductions on the test set for the heuristic, agents $A_5$ and agent $A_6$.}
	\label{fig: daily based bill improvement}
\end{figure}

We show in Fig. \ref{fig: daily based bill improvement} and Table \ref{tab:DDM_results} the results for problem (DDM) of the heuristic and two RL agents: $A_5$ with state space 
\begin{equation} \label{eq:s5}
    S_4=\mathcal{T} \times \mathcal{S}o\mathcal{C}_d \times \Delta_d \times \mathcal{P}_d
\end{equation}
and completed by the lazy controller during the tests, and $A_6$ the same with $A_5$ but completed by the heuristic controller during the tests.
Each DPI loop has a computational complexity of $O(\left(N_{train} |S_4|\right)$.

We can see that the results of agents $A_5$ and $A_6$ have very small differences, and they both reduces the bill compared to the heuristic controller, which failed at least once to reduce the total bill.

\begin{figure}
    \centering
     \hspace*{-0.3cm}
     \begin{subfigure}[b]{0.5\columnwidth}
         \centering
         \includegraphics[width=\textwidth]{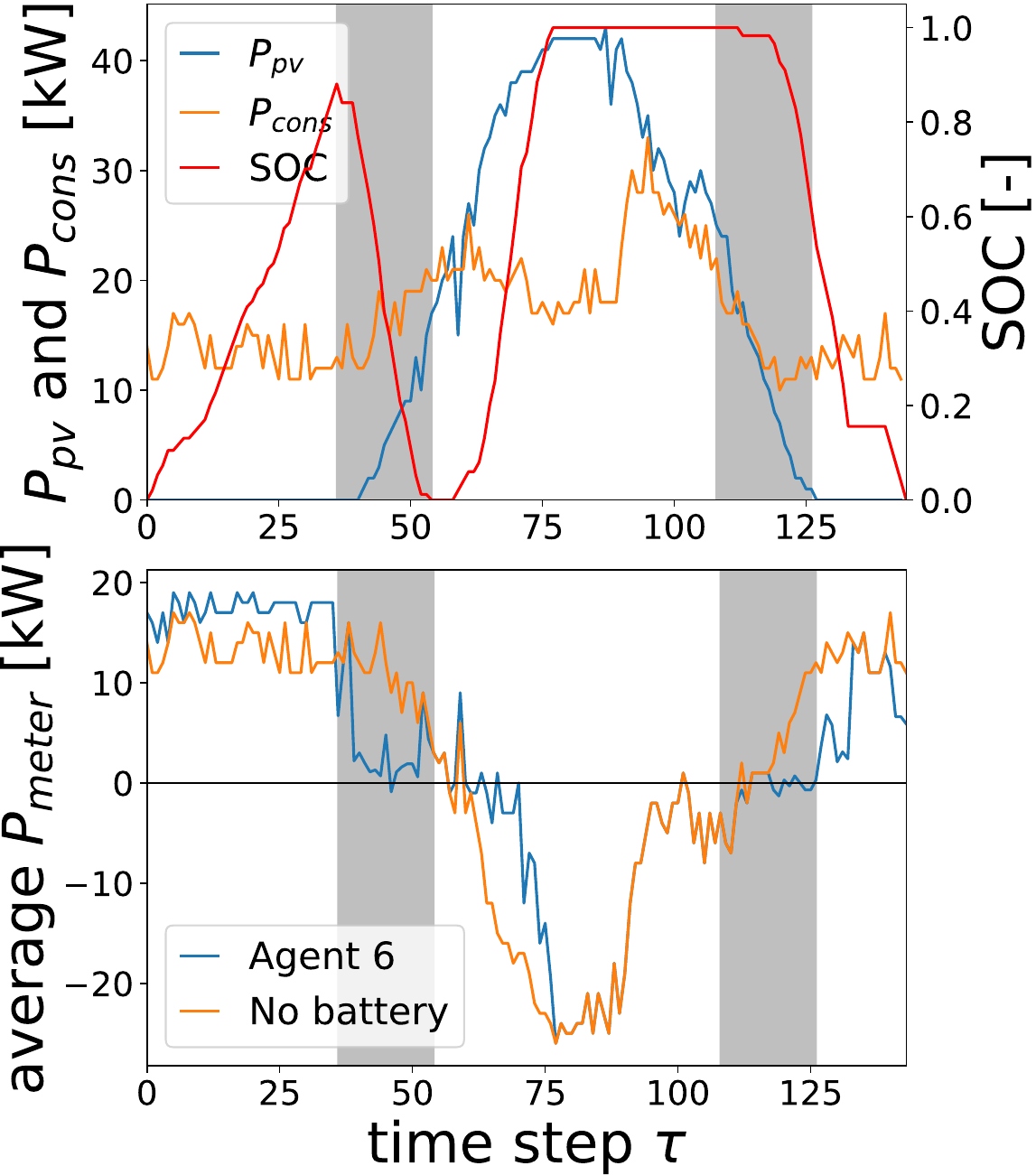}
         \caption{Day with small bill reduction}
         \label{fig: bad day daily DC}
     \end{subfigure}
     \begin{subfigure}[b]{0.5\columnwidth}
         \centering
         \includegraphics[width=\textwidth]{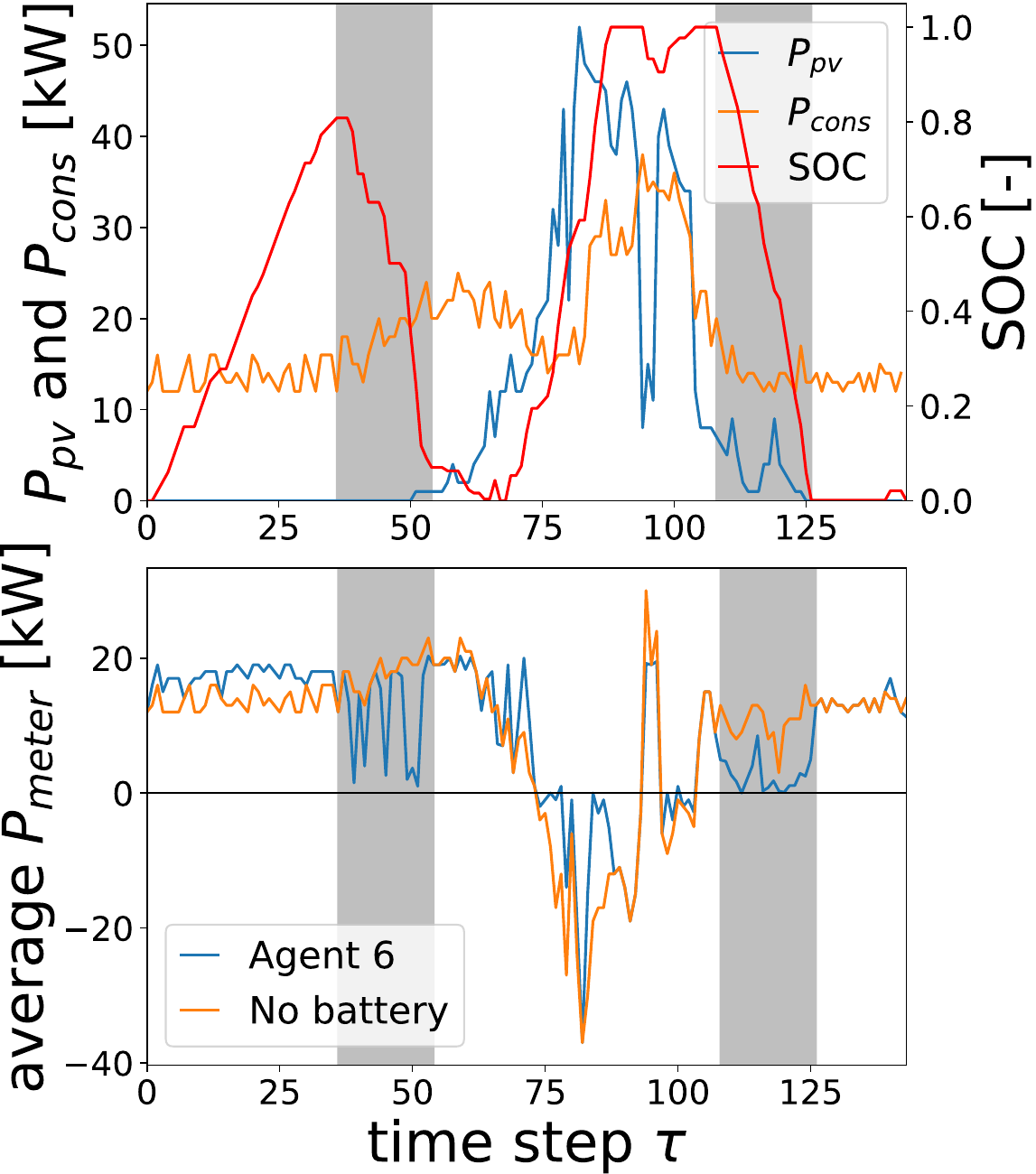}
         \caption{Day with big bill reduction}
         \label{fig: good day daily DC}
    \end{subfigure}
    \caption{Days with small and big bill reduction. Top figures show the power generation, the electricity demand and the battery charge. Bottom figures represent $P_{meter}$ with and without battery. Gray areas are on-peak hours.}
    \label{fig: good and bad day for daily DC}
\end{figure}

We illustrate typical days with big and small bill reductions in Fig. \ref{fig: good and bad day for daily DC} and show detailed bill reduction values in Table \ref{tab:DDM_results}. 
For the day with a big bill reduction in Fig. \ref{fig: good day daily DC}, the production falls suddenly during high consumption period (between time steps 75 and 105), and leads to a high peak demand. 
The battery charged during low consumption period (at the beginning of the day) then allows to significantly reduce the DC.
For the day with a small bill reduction in Fig. \ref{fig: bad day daily DC}, the battery is charged before the first on-peak hours to absorb the on-peak consumption and to deal with potential important peak demands. 
This strategy appears natural since the peak demand during high consumption period risks to increase significantly the electricity bill (like in the big bill reduction case). Unfortunately, as shown by Table \ref{tab:DDM_results}, this decision may increase slightly the DC when the high loads are totally covered by solar production.

\begin{remark}
We notice that the battery is not fully charged before the first on-peak hour. 
The speed of charge is limited to keep a small demand peak. 
This indicates that the benefits of having a bigger battery may be limited. 
In addition, this indicates that the problem decomposition described in section \ref{problem decomposition} may be improved by recoupling the separate days. We are to discuss this in section \ref{recoupling subsection}.
\end{remark}

 \begin{table}[]
    \centering
    \begin{tabular}{|c|c|c|c|}
    \hline
   Reduction & bill ($\$$) & energy charge ($\$$) & DC ($\$$) \\
   \hline
   Average agent $A_5$ &  3.26 & 2.01 & 1.25\\
   Average agent $A_6$ &  3.34 & 2.02 & 1.32\\
   Day with big reduction & 7.64 & 2.64 & 5.00\\
   Day with small reduction & 0.99 & 1.99 & -1.0\\
   \hline
    \end{tabular}
    \caption{Bill reduction results for Problem (DDM)}
    \label{tab:DDM_results}
\end{table}

\subsection{Results of problem (MDM)}
\label{subsection: monthly DC}
\begin{figure}
    \centering
    \includegraphics[width=0.8\columnwidth]{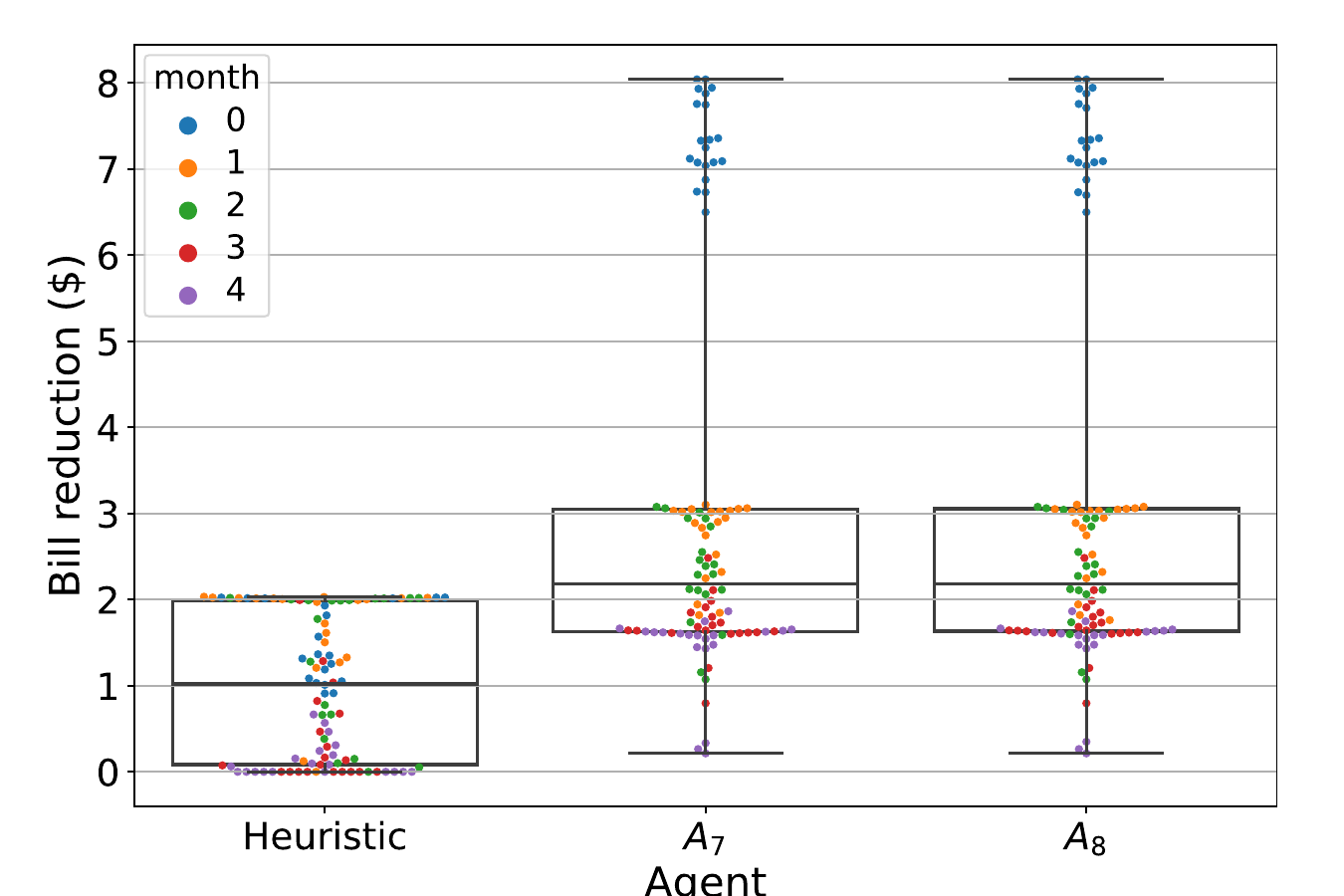}
    \caption{Swarm plot representing the daily bill reductions on the test set for the heuristic, agent $A_7$ and agent $A_8$.}
	\label{fig: month-based bill minimization}
\end{figure}

Similar to the previous section, we consider two RL agents:
agent $A_7$ lives in $S_5$ defined by \eqref{eq:s5} completed by the lazy controller, 
and agent $A_8$ completed by the heuristic controller.
The results for problem (MDM) are shown in Fig. \ref{fig: month-based bill minimization} and Table \ref{tab:MDM_results}. 
To be able to compare results of problem (DDM) and problem (MDM), the daily bill for the problem (MDM) is defined as the sum of the energy charge and the MDC divided by 20 (the number of working days).
Compared to the results of problem (DEM), maximum gains are again much higher. 
Again, we observe that our agents are systematically beneficial, and the difference between the use or not of the heuristic remains marginal.

 \begin{table}[]
    \centering
    \begin{tabular}{|c|c|c|c|}
    \hline
   Reduction & bill ($\$$) & energy charge ($\$$) & DC ($\$$) \\
   \hline
   Daily average agent $A_7$ &  3.09 & 2.29 & 0.80\\
   Daily average agent $A_8$ &  3.10 & 2.30 & 0.80\\
   Day with big reduction & 8.04 & 3.04 & 5.00\\
   Day with small reduction & 0.21 & 0.71 & 0.50\\
   Month 0 & 7.34 & 2.34 & 5.00 \\
   Month 1 & 2.73 & 2.73 & 0.00 \\
   Month 2 & 2.34 & 2.34 & 0.00\\
   Month 3 & 1.69 & 2.19 & -0.50 \\
   Month 4 & 1.41 & 1.91 & -0.50 \\
   \hline
    \end{tabular}
    \caption{Bill reduction results for Problem (MDM)}
    \label{tab:MDM_results}
\end{table}

\begin{figure}
    \centering
    
     \hspace*{-0.3cm}
     \begin{subfigure}[b]{0.5\columnwidth}
         \centering
         \includegraphics[width=\textwidth]{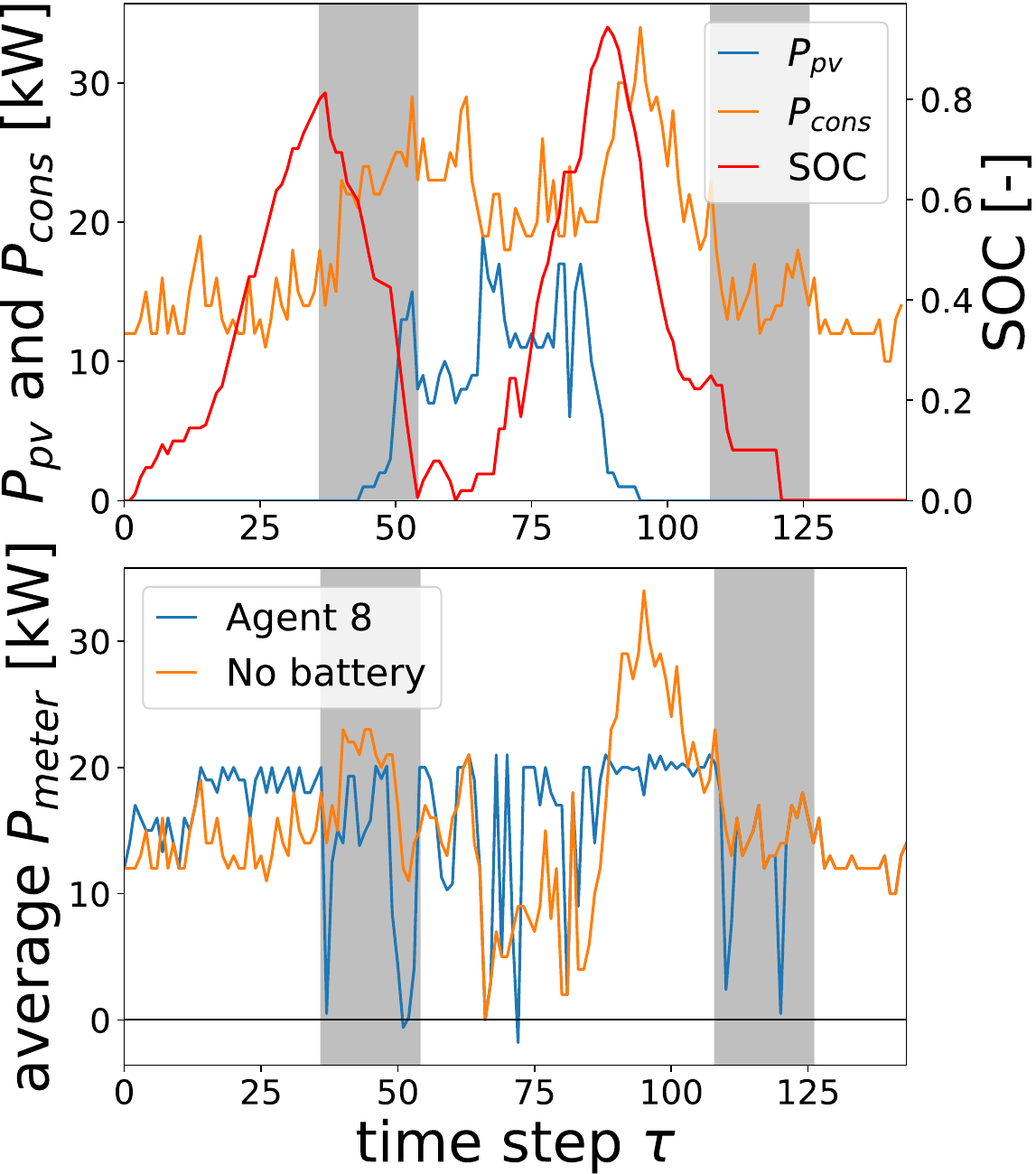}
         \caption{Day with small bill reduction}
         \label{fig: bad day monthly DC}
     \end{subfigure}
     \begin{subfigure}[b]{0.5\columnwidth}
         \centering
         \includegraphics[width=\textwidth]{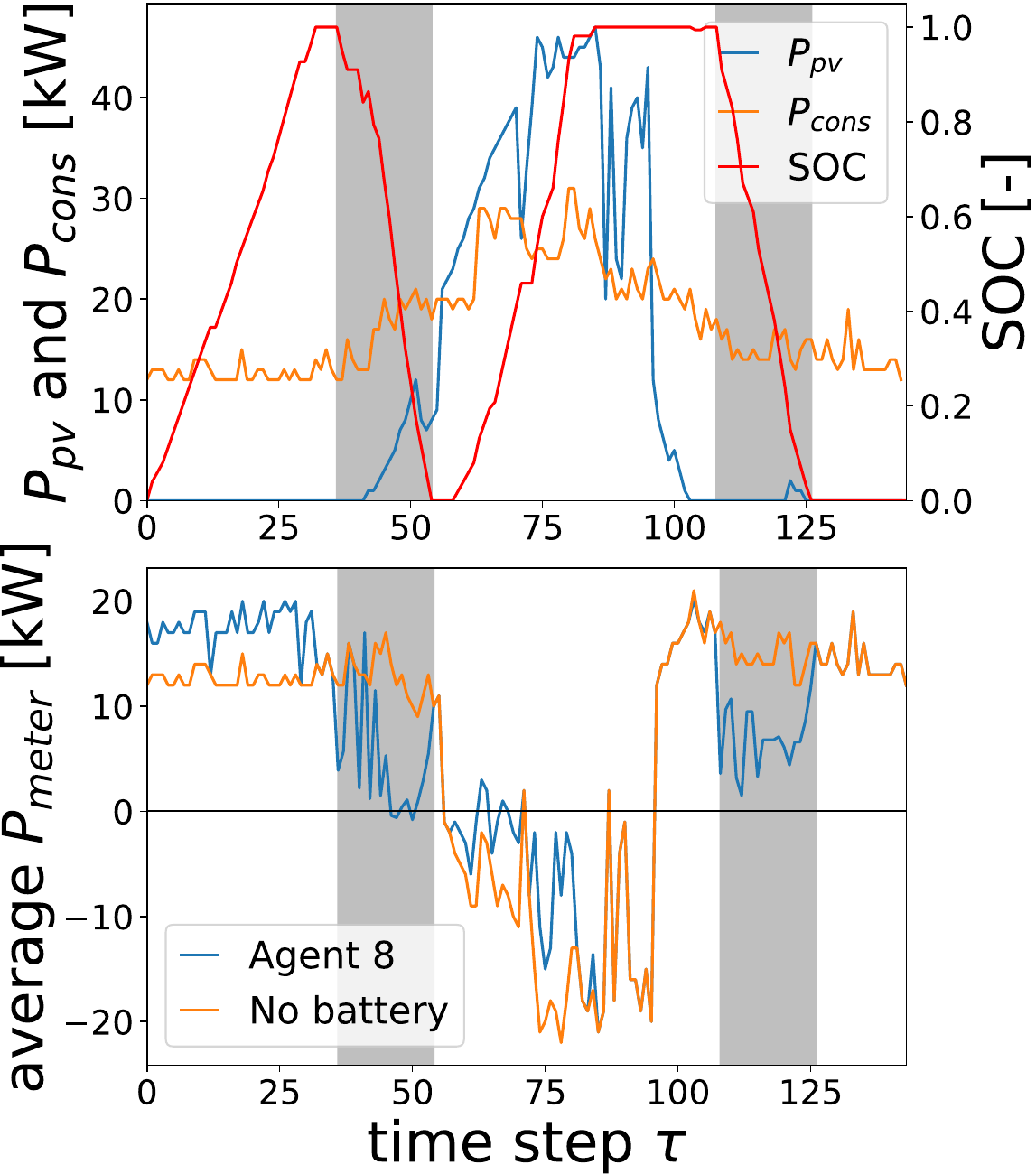}
         \caption{Day with big bill reduction}
         \label{fig: good day monthly DC}
    \end{subfigure}
    \caption{Days with small and big bill reduction. Top figures show the power generation, the electricity demand and the battery charge. Bottom figures represent $P_{meter}$ with and without battery. Gray areas are on-peak hours.}
    \label{fig: good and bay day for monthly DC}
\end{figure}

We show in Fig. \ref{fig: good and bay day for monthly DC} typical days with small and big bill reductions. 
The day with a small bill reduction shown in Fig. \ref{fig: bad day monthly DC} is the first day of the month. 
We can see that the agent's behavior is cautious: the power demand never exceeds $20 kW$. The lost on the DC is caused by another day of the month. 
On contrary, the behavior of the agent is much less conservative for the day with big bill reduction shown in Fig. \ref{fig: good day monthly DC}. 
It is a day in the middle of the month, and a higher peak of $20 kW$ has been measured earlier in this month.
We can see that before the first on-peak hours, the net power demand is instantly bigger and allows the battery to be fully charged. This energy storage leads then to smaller net power demands during the first on-peak hours. 

These results correspond to the intuition. Before a high peak is reached during the BP, the conservative behavior limits the economies with energy charges. Afterwards, electricity demands are less restrictive and the peak load to the grid is less reduced.

\subsection{Comparison of problems (DDM) and (MDM)}
\begin{figure}
    \centering
    \includegraphics[width=0.8\columnwidth]{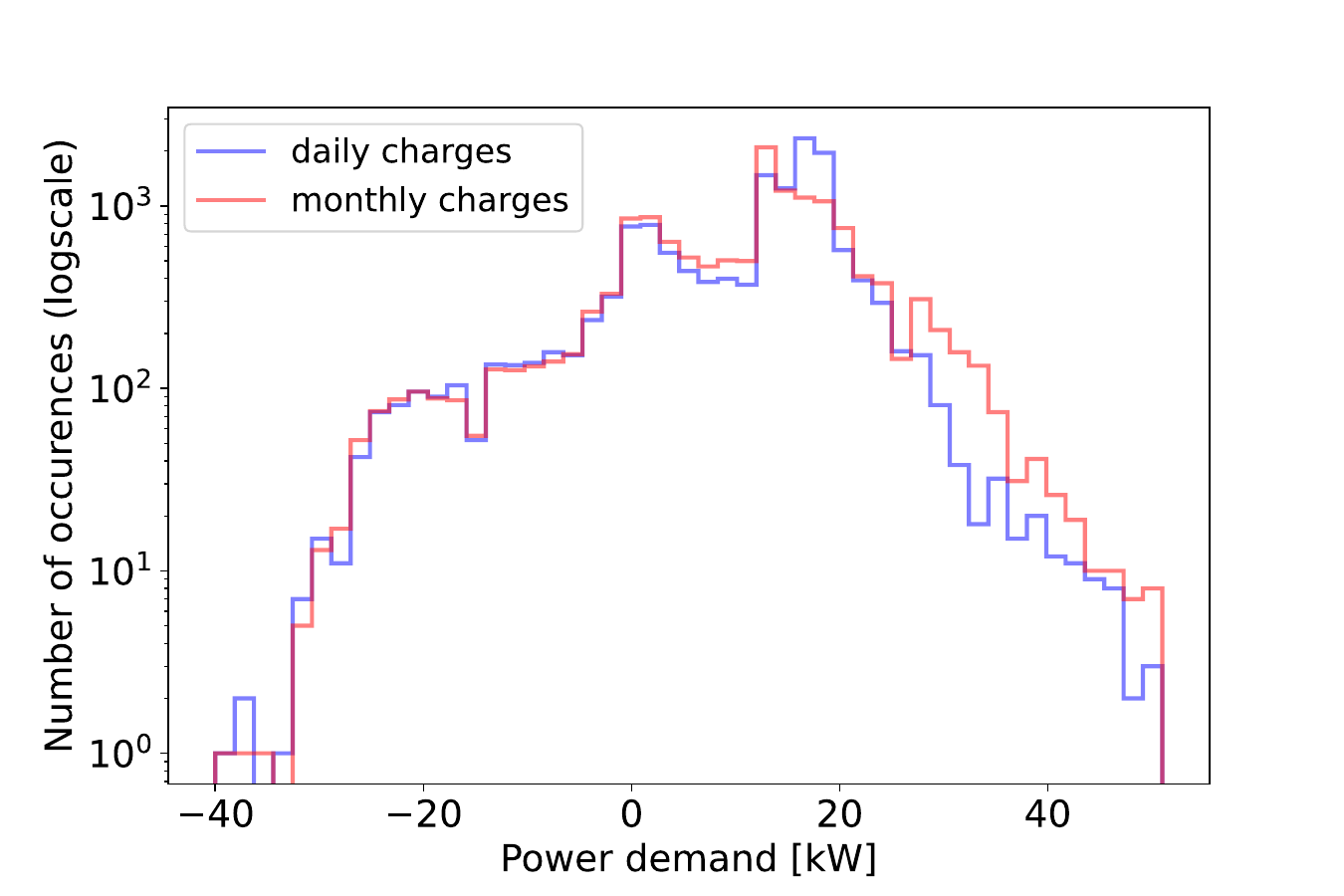}
    \caption{Histograms of the power provided by the grid (kW) measured on intervals of 10 minutes for MDC and DDC.}
	\label{fig: BP comparison}
\end{figure}

From Table \ref{tab:DDM_results} and Table \ref{tab:MDM_results}, we can see that the average daily bill reduction is slightly bigger with problem (DDM).
Indeed, for problem (MDM), the test month $0$ is monopolizing the very high bill reductions, whereas the bill reductions remain very small for other months.
In Fig. \ref{fig: BP comparison}, we show the distribution of the measured net power demands $P_{meter}$. 
We can observe that for problem (MDM), the numbers of high peaks are larger than those for problem (DDM). 
These results suggest that using the DDC is more likely to decrease the load peaks on the grid.

A simple explanation is that for a given peak, the agent is not penalized as long as he does not exceed it. As explained in the previous section, once a high peak has been reached in case of MDC, the consumer becomes less motivated to reduce its demand peak for the rest of the month, and thus may lead to higher peaks.

\subsection{Re-coupling}
\label{recoupling subsection}
There may be room to increase the gains by re-coupling separate days. 
A first attempt is to charge the battery at the end of the day during the last off-peak hours (before the next on-peak hours), while being careful not to exceed the peak demand reached during the day.
By applying this recoupling strategy on our trained agents, we observe an additional bill reduction (of $\$3.65$) for problem (DDM).
However, the performance degrades for problem (MDM) (with a bill reduction of $\$3.10$). 
One reason for this phenomenon may be the use of higher charging powers for problem (MDM).
Indeed, for the nonlinear battery model \eqref{eqn:x_dot}, the charging speed $\dot{x}$ grows sub-linearly when the charging power $P_{c,bat}$ increases and decreases sup-linearly when $P_{d,bat}$ increases. Consequently, the faster the battery is charged or discharged, the lower the efficiency of the battery's power transition, and therefore the greater the power loss.

\section{CONCLUSIONS AND FUTURE WORKS}
\label{section 6}
We designed a new RL algorithm for energy and demand charges minimization. As no hyperparameters are involved and no production or consumption forecasts are required, this method is particularly simple to implement.
We showed that significant benefits can be achieved 
through intelligent battery steering, especially when demand charges are involved. 
We compared DDCs to MDCs and showed that the former can lead to lower power demand peaks and offers bigger potential economic gains.

The recoupling of days remains an open question that has been highlighted by the consideration of battery non-linearities. Future work could focus on this aspect, as well as on implementation of RL methods for continuous state and action spaces.

\section{ACKNOWLEDGEMENTS}
This work was supported by the Direction Générale de l’Armement.
\Urlmuskip=0mu plus 1mu\relax
\bibliographystyle{IEEEtran}
\bibliography{main}

\begin{thebibliography}{10}
\providecommand{\url}[1]{#1}
\csname url@samestyle\endcsname
\providecommand{\newblock}{\relax}
\providecommand{\bibinfo}[2]{#2}
\providecommand{\BIBentrySTDinterwordspacing}{\spaceskip=0pt\relax}
\providecommand{\BIBentryALTinterwordstretchfactor}{4}
\providecommand{\BIBentryALTinterwordspacing}{\spaceskip=\fontdimen2\font plus
\BIBentryALTinterwordstretchfactor\fontdimen3\font minus
  \fontdimen4\font\relax}
\providecommand{\BIBforeignlanguage}[2]{{%
\expandafter\ifx\csname l@#1\endcsname\relax
\typeout{** WARNING: IEEEtran.bst: No hyphenation pattern has been}%
\typeout{** loaded for the language `#1'. Using the pattern for}%
\typeout{** the default language instead.}%
\else
\language=\csname l@#1\endcsname
\fi
#2}}
\providecommand{\BIBdecl}{\relax}
\BIBdecl

\bibitem{Basics_Electricity_Price_Formation}
P.~Morthost, S.~Ray, J.~Munksgaard, and A.~Sinner, ``Wind energy and
  electricity prices. exploring the 'merit order effect','' 2010.

\bibitem{Relation_price_pollution}
F.~A. Felder, ``Examining electricity price suppression due to renewable
  resources and other grid investments,'' \emph{The Electricity Journal},
  vol.~24, no.~4, pp. 34--46, 2011.

\bibitem{Demand_Side_Management}
C.~W. Gellings, ``The concept of demand-side management for electric
  utilities,'' \emph{Proceedings of the IEEE}, vol.~73, no.~10, pp. 1468--1470,
  1985.

\bibitem{Time_of_Use_pricing}
G.~R. Newsham and B.~G. Bowker, ``The effect of utility time-varying pricing
  and load control strategies on residential summer peak electricity use: A
  review,'' \emph{Energy policy}, vol.~38, no.~7, pp. 3289--3296, 2010.

\bibitem{Demand_Charges}
J.~A. McLaren, P.~J. Gagnon, and S.~Mullendore, ``Identifying potential markets
  for behind-the-meter battery energy storage: A survey of us demand charges,''
  National Renewable Energy Lab.(NREL), Golden, CO (United States), Tech. Rep.,
  2017.

\bibitem{PGOptionStorage}
\BIBentryALTinterwordspacing
``{PG}\&{E}’s {Option} {S} (for {Storage}) {Daily} {Demand} {Charge} {Rate}
  – {Here}'s what you should know.'' [Online]. Available:
  \url{https://www.energytoolbase.com/newsroom/blog/pges-option-s-for-storage-daily-demand-charge-rate-heres-what-you-should-know}
  (visited on Sep. 5, 2023).
\BIBentrySTDinterwordspacing

\bibitem{Augmented_state_space}
M.~J. Risbeck and J.~B. Rawlings, ``Economic model predictive control for
  time-varying cost and peak demand charge optimization,'' \emph{IEEE
  Transactions on Automatic Control}, vol.~65, no.~7, pp. 2957--2968, 2019.

\bibitem{Radiation_forecasts}
G.~Van~Kriekinge, C.~De~Cauwer, N.~Sapountzoglou, T.~Coosemans, and
  M.~Messagie, ``Peak shaving and cost minimization using model predictive
  control for uni-and bi-directional charging of electric vehicles,''
  \emph{Energy reports}, vol.~7, pp. 8760--8771, 2021.

\bibitem{Energy_forecasts}
C.~Croonenbroeck and G.~Stadtmann, ``Renewable generation forecast
  studies--review and good practice guidance,'' \emph{Renewable and Sustainable
  Energy Reviews}, vol. 108, pp. 312--322, 2019.

\bibitem{Demand_forecasts}
Y.~Wang, D.~Gan, M.~Sun, N.~Zhang, Z.~Lu, and C.~Kang, ``Probabilistic
  individual load forecasting using pinball loss guided lstm,'' \emph{Applied
  Energy}, vol. 235, pp. 10--20, 2019.

\bibitem{xiang2015robust}
Y.~Xiang, J.~Liu, and Y.~Liu, ``Robust energy management of microgrid with
  uncertain renewable generation and load,'' \emph{IEEE Transactions on Smart
  Grid}, vol.~7, no.~2, pp. 1034--1043, 2015.

\bibitem{mo2021optimal}
Y.~Mo, Q.~Lin, M.~Chen, and S.-Z.~J. Qin, ``Optimal online algorithms for
  peak-demand reduction maximization with energy storage,'' in
  \emph{Proceedings of the twelfth ACM international conference on future
  energy systems}, 2021, pp. 73--83.

\bibitem{EC1}
V.~Fran{\c{c}}ois-Lavet, D.~Taralla, D.~Ernst, and R.~Fonteneau, ``Deep
  reinforcement learning solutions for energy microgrids management,'' in
  \emph{European Workshop on Reinforcement Learning (EWRL 2016)}, 2016.

\bibitem{EC2}
J.~Cao, D.~Harrold, Z.~Fan, T.~Morstyn, D.~Healey, and K.~Li, ``Deep
  reinforcement learning-based energy storage arbitrage with accurate
  lithium-ion battery degradation model,'' \emph{IEEE Transactions on Smart
  Grid}, vol.~11, no.~5, pp. 4513--4521, 2020.

\bibitem{EC3}
C.~Guan, Y.~Wang, X.~Lin, S.~Nazarian, and M.~Pedram, ``Reinforcement
  learning-based control of residential energy storage systems for electric
  bill minimization,'' in \emph{2015 12th Annual IEEE Consumer Communications
  and Networking Conference (CCNC)}.\hskip 1em plus 0.5em minus 0.4em\relax
  IEEE, 2015, pp. 637--642.

\bibitem{Policy_Iteration}
R.~S. Sutton and A.~G. Barto, \emph{Reinforcement learning: An
  introduction}.\hskip 1em plus 0.5em minus 0.4em\relax MIT press, 2018.

\bibitem{Dynamic_Battery}
L.~Guzzella, A.~Sciarretta, L.~Guzzella, and A.~Sciarretta, ``Vehicle energy
  and fuel consumption--basic concepts,'' \emph{Vehicle propulsion systems:
  introduction to modeling and optimization}, pp. 13--46, 2013.

\end{thebibliography}

\end{document}